\documentclass[11pt, twoside, a4paper]{article}
\usepackage[utf8]{inputenc}
\usepackage[english]{babel}

\usepackage{amssymb,amsmath,amsthm}
\usepackage{mathtools} 
\usepackage{bigints}
\usepackage{breakcites}

\usepackage{dsfont}
\usepackage[mathscr]{eucal} % mathcal style
\usepackage{graphicx}

\usepackage{xcolor}
\usepackage[hidelinks]{hyperref}
\hypersetup{
 colorlinks,
 linkcolor={red!50!black},
 citecolor={blue!50!black},
 urlcolor={blue!80!black}
}
\usepackage{pdfsync}
\usepackage{enumerate} % Fancy enumeration lists 
\usepackage{esint}

\theoremstyle{plain}
\newtheorem{theorem}{Theorem}[section]
\newtheorem{definition}[theorem]{Definition}

\newtheorem{proposition}[theorem]{Proposition}
\theoremstyle{remark}
\newtheorem{remark}[theorem]{Remark}

\numberwithin{equation}{section}
\allowdisplaybreaks

\def\R{{\mathbb R}}% real numbers
\newcommand{\dd}{\,\mathrm{d}}

\renewcommand{\div}{\operatorname{div}}

\renewcommand{\leq}{\leqslant}
\renewcommand{\geq}{\geqslant}
\renewcommand{\epsilon}{\varepsilon}

\newcommand{\dt}{\partial_t}
\newcommand{\eps}{\varepsilon}
\newcommand{\bigO}{\mathcal{O}}
\newcommand{\ap}{\alpha_+}
\newcommand{\am}{\alpha_-}
\newcommand{\rp}{\rho_+}
\newcommand{\rmi}{\rho_-}
\newcommand{\gp}{\gamma_{+}}
\newcommand{\gm}{\gamma_-}

\title{Low-Mach-number limit for two-phase flows} 
\author{Cassandre Lebot}
\def\adrese{
\noindent CNRS, LAMA, ISTerre, Univ. Savoie Mont Blanc, 73000 Chamb\'ery, France.\\
{\it Email address:}
\texttt{Cassandre.Lebot@univ-smb.fr}
}
\date{\today}

\begin{document}
\maketitle

\begin{abstract}
This paper is devoted to the formal study of the low-Mach-number limit for solutions of the compressible Navier-Stokes or Euler equations for different types of fluids. We first review the different results obtained in the case of flows consisting of one phase. Then, we focus on the low-Mach-number limit for two-phase flows, considering different types of systems: with an algebraic closure or a PDE closure for the pressure, with one single or two different velocities, without or with entropy.
\end{abstract}

\tableofcontents

\section{Introduction}

The aim of this proceeding is to introduce the concept of the low-Mach-number limit for compressible multiphase flows at the formal level. The Mach number is a dimensionless quantity defined as the ratio of the flow velocity to the speed of sound in the fluid. When the Mach number is low, compressibility effects are negligible and the flow may be approximated as incompressible. In the past, several authors have worked on the passage from compressible equations to an incompressible system when the Mach number tends to zero. The motivations are several. First, we better understand incompressible equations and we may use the passage to the limit and comparison to understand more properties on compressible equations. Moreover, low-Mach-number equations have several advantages in numerical analysis: simplification of the algorithms, smaller time step for the numerical integration. For more details on the numerical challenges for low-Mach-number systems, we refer the reader to \cite{battisti2025linearly}.

However, this question of the zero-Mach-number limit for multiphase flows is still to be examined. To the best of our knowledge, \cite{varsakelis2011low} and \cite{narbona2024two} are two of the very few articles that studied this question. We may also refer to \cite{fanelli2025incompressible} for related incompressible limits in a reduced compressible MHD system. The present contribution aims at giving an overview of this subject by discussing the formal derivation of the zero-Mach-number limit for several types of multiphase flow models. A rigorous mathematical justification of one of these limits is provided in \cite{lebot2026low}.

This paper is organized as follows. In Section~\ref{seconeP}, we give a brief review of the low-Mach-number limit for one-phase flow. In Section~\ref{sec2P1U}, we focus on two-phase flows system with the same velocity for both phases and derive the low-Mach-number equations for two-phase systems considering two types of pressure closure: an algebraic closure and a PDE closure. We discuss first the isentropic case and then the non-isentropic case. In Section~\ref{sec2P2U}, we derive the low-Mach-number equations for two-phase flows system with two distinct velocities. Finally, Section~\ref{Opb} contains a conclusion, some open problems and remarks about topics that we have not considered in this paper.

\section{Low-Mach number limit for one-phase flow}\label{seconeP}

We give here a very short review of some results on the low Mach number limit for compressible flows. We refer the interested reader for more detailed reviews to \cite{schochet2005mathematical}, \cite{gallagher2005resultats} or \cite{alazard2008minicourse}.

\subsection{The isentropic case}

In this part, we restrict ourselves to the case of isentropic fluids, caracterized by a constant entropy, and we assume that the flow evolves in the domain $\Omega$ from the $d$-dimensional space $\R^d$, with $d \geq 2$. We focus on domains that are either periodic or bounded. In bounded domains, we assume the Dirichlet boundary conditions, namely $u=0$ on $\partial \Omega$.

We write $t\in \R^+$ for the time variable and $x \in \Omega$ denotes the space variable. We consider the following governing equations of a compressible flow: 
\begin{equation}
\left\{
\begin{aligned}\label{Sys1}
    &\dt \rho + \div (\rho u) = 0,\\
    &\dt (\rho u) + \div (\rho u \otimes u) + \nabla p(\rho) - \mu \Delta u - (\lambda + \mu)\nabla \div u =0
\end{aligned}
\right.
\end{equation}
with initial conditions
\begin{align*}
    \rho_{\lvert t=0} = \rho_0 \geq 0 \quad \text{and} \quad u_{\lvert t=0} = u_0
\end{align*}
and where $\rho = \rho(t,x) \geq 0$ represents the mass density, $u = u(t,x)$ the velocity vector, $p=p(\rho(t,x))$ the pressure. The viscosity coefficients $\lambda$ and $\mu$ satisfy $\mu >0$ and $\lambda + 2\mu>0$ and are assumed to be constant (this can be relaxed). We consider a barotropic fluid with a $\gamma$ pressure law: 
\[
p(\rho) = \rho^\gamma\text{, with }\gamma >1.
\]

There are several ways on how to introduce the Mach number in the compressible flow equations. First, we can consider a varying equations of state and also varying transport coefficients, like viscosity coefficients for viscous fluids (see \cite{ebin1977motion}, \cite{schochet1986compressible}). In that approach, a particular fluid is viewed as ``nearly incompressible''. 

The alternative way to derive the compressible Navier-Stokes equations in the small Mach number regime is to consider flows on varying time or spatial scale. This amounts to use a dimensional analysis, by rewriting the equations in their dimensionless form. One way to do it is by scaling each variable by its characteristic value (see \cite{klein2001asymptotic}). The other way is to use a change of independent variables (see \cite{alazard2006low}, \cite{lions1998incompressible}). This is the approach we take in the following development (for more details, see the introduction in \cite{lions1996mathematical}).\\

We denote by $\eps$ the Mach number. In contrast with the usual physical scaling, where the Mach number appears explicitly in front of the pressure and viscous terms, we adopt a rescaling of time and space variables. This approach is asymptotically equivalent and corresponds to the same low-Mach-number regime. It allows us to avoid artificial $\eps$-dependent viscosity coefficients and simplifies the formal asymptotic expansion, while preserving the correct balance between inertial, viscous, and pressure forces.

The incompressible limit corresponds to flow regimes in which the speed of the fluid becomes negligible compared to the speed of sound. In order to capture this low-Mach-number regime, we rescale the velocity by writing $u = \epsilon u^\epsilon$, where $u^\epsilon$ is of order one. Since the velocity is small, fluid particles move only over long time scales. More precisely, a particle travels a distance of order $1$ over a time scale of order $1/\epsilon$, and a distance of order $1/\epsilon$ over a time scale of order $1/\epsilon^2$, meaning that the relevant time scale is $1/\eps^2$. This shows that the relevant dynamics take place on large spatial and long temporal scales. For this reason, we introduce the rescaled variables
\[
t^\epsilon = \epsilon^2 t, \qquad x^\epsilon = \epsilon x,
\]
so that both $t^\epsilon$ and $x^\epsilon$ are of order one (see \cite{danchin2005low}). Thus, we scale $u$ and $\rho$ in the following way
\begin{align}
    u(t,x) = \eps u^\eps(t^\eps,x^\eps), \quad  \rho(t,x) = \rho^\eps(t^\eps,x^\eps).\label{changevar}
\end{align}

We can rewrite System \eqref{Sys1}: 
\begin{equation}
\left\{
\begin{aligned}\label{Sys2}
    &\partial_{t^\eps} \rho^\eps + \div (\rho^\eps u^\eps) = 0,\\
    &\partial_{t^\eps} (\rho^\eps u^\eps) + \div (\rho^\eps u^\eps \otimes u^\eps) + \frac{1}{\eps^2}\nabla p(\rho^\eps) - \mu^\eps \Delta u^\eps - (\lambda^\eps + \mu^\eps)\nabla \div u^\eps =0.
\end{aligned}
\right.
\end{equation}

The initial data are $\rho^\eps_{\lvert t=0 }=\rho^\eps_0$ and $u^\eps_{\lvert t=0 }=u^\eps_0$. We will see that it can be relevant to use initial data that depend on $\eps$.

Formally, when we take the limit $\eps \to 0$ in System \eqref{Sys2}, we cancel the terms with the highest factor in $\eps$ in the equations. Therefore from the momentum equation we obtain when $\eps$ goes to zero
\begin{align}
    \nabla p(\rho^\eps) \to 0 \label{cond1}.
\end{align}

Let's call $C$ a constant in space such that the limit density $\rho^0$ satisfies, because of \eqref{cond1}:
\begin{align*}
    p(\rho^0) = C(t),\\
    \rho^0 = C(t)^{1/\gamma}.
\end{align*}

From the mass equation we get
\begin{align*}
    \frac{1}{\gamma} \frac{C'(t)}{C(t)} + \div u^0 = 0.
\end{align*}

We now integrate in space over $\Omega$, that we have assumed periodic or bounded with the Dirichlet boundary conditions:
\begin{align*}
    \int_\Omega \frac{1}{\gamma} \frac{C'(t)}{C(t)} = 0.
\end{align*}
We deduce that $C'(t)=0$ and thus that $\rho^\eps$ must converge to a positive constant. Note that in the case of the whole space, this depends on the condition that we use on $\rho^0$ at infinity.

Now passing to the limit in the mass equation we obtain
\begin{align}
    \div u^\eps \to 0\label{cond2}.
\end{align}

Finally, coming back to the momentum equation, we find that the limit velocity $u^0$ must satisfy the incompressible Navier-Stokes equations
\begin{equation*}
\left\{
\begin{aligned}
   &\div u^0 =0,\\
   &\dt (\rho^0 u^0) + \div (\rho^0 u^0 \otimes u^0) - \mu \Delta u^0 + \nabla \pi^0 = 0, 
\end{aligned}
\right.
\end{equation*}
with $\pi^0$ being the ``limit'' of $\frac{(\rho^\eps)^\gamma-(\rho^0)^\gamma}{\eps^2}$.

\begin{remark}
We can also rewrite the momentum equation as
$$ \rho^0 \Big( \dt u^0 + (u^0 \cdot \nabla) u^0 \Big)- \mu \Delta u^0 + \nabla \pi^0 = 0.  $$
\end{remark}

We end up with a compatibility condition on the initial data: it is more or less difficult to pass to the limit in System \eqref{Sys2} if the initial data verifies itself the conditions \eqref{cond1}-\eqref{cond2} or not. Clearly, it is easier to pass to the limit if it does.  Hence, we introduce the notion of ``well--prepared'' initial data. 

\begin{definition}\label{defwpd} We say that the initial data $(u_0^\eps, \rho_0^\eps)$ are well-prepared if they satisfy
\begin{align*}
    \div u_0^\eps = \bigO(\eps^2), \quad  \rho_0^\eps = \bar{\rho_0} + \bigO(\eps^2)
\end{align*}
with $\bar{\rho_0}$ a constant.
\end{definition}

Otherwise, we say that the initial data are ``ill-prepared'' or we call them ``general'' initial data. In that case, it actually means that sounds waves are generated initially and propagate with high speed in the domain (speed of order $1/\eps$). In that paper, we will not consider this case and focus only on problems with well-prepared initial data. Under this assumption, the nature of the domain, periodic or bounded, does not affect the leading-order asymptotics, and the formal derivation of the low-Mach-number limit proceeds in the same way.

Historically, the study of the low-Mach number limit in the compressible Navier-Stokes equations began with the case of well-prepared data. The first general proof of the transition from the equations of compressible fluids to the equations of incompressible fluids comes from \cite{ebin1977motion}. Then \cite{klainerman1981singular}, \cite{klainerman1982compressible} also studied the existence time of the solutions of System \eqref{Sys2} (its uniformity in $\eps$) and the passing to the limit, using symmetric hyperbolic system theory in the whole space and in the periodic case respectively. We can also note the work of \cite{schochet1986compressible} in a bounded domain with periodic boundary conditions that studies the convergence using filtering techniques.

Later, more works have dropped the assumption of well-prepared data. Considering the whole space, \cite{ukai1986incompressible} gives the first proof of the transition from compressible Navier-Stokes equations to incompressible Navier-Stokes equations when the initial data is more general. However the uniform convergence does not hold in $t=0$ because of the development of an initial layer (see also \cite{iguchi1997incompressible}, \cite{isozaki1986singular}). With periodic boundary conditions, several authors have studied the convergence using the filtering method (\cite{danchin2002zero}, \cite{grenier1997oscillatory}, \cite{gallagher1998asymptotic}).

P.-L. Lions, B. Desjardins, E. Grenier and N. Masmoudi also show in several works that we can pass to the limit when $\eps$ goes to zero using uniform estimates coming from the energy inequality. We remind the reader the energy inequality for Navier-Stokes equation, as it will be interesting to compare it later to the ones of two-phase flows systems:
\begin{multline} \label{energie1}
    \int_\Omega \biggl[\frac{1}{2} \rho^\eps \lvert u^\eps \rvert^2(t,\cdot) + \frac{1}{\eps^2 (\gamma-1)} (\rho^\eps)^\gamma (t,\cdot)\biggr] + \int_0^t\!\!\!  \int_\Omega \biggl[\mu^\eps \lvert \nabla u^\eps \rvert^2 + \lambda^\eps (\div u^\eps)^2\biggr] \\
    \leq \int_\Omega \biggl[\frac{1}{2} \rho_0^\eps \lvert u_0^\eps \rvert^2 + \frac{1}{\eps^2(\gamma-1)} (\rho_0^\eps)^\gamma\biggr].
\end{multline}
We refer to Subsection~\ref{subE} for more details on the computations of this energy.

They give a proof in \cite{lions1998incompressible} for the case of periodic boundary conditions, and in \cite{desjardins1999low} for weak solutions for isentropic fluids in the whole space. In the latter, they improve the convergence result using Strichartz' estimates. For the case of a bounded domain with homogeneous Dirichlet conditions, \cite{desjardins1999incompressible} show that the convergence is strong because of acoustic waves absorbed in boundary (using a geometric assumption on the domain). Let us note also the work of \cite{secchi2000singular} that shows the weak convergence to solutions of incompressible flows using the antisymmetry of the acoustic operator when the condition $u \cdot \mathbf{n} =0$ is set in the boundary (and the convergence is shown to be strong in the case of well-prepared data). For a new approach to the low-Mach-number limit using dissipative measure-valued solutions, we refer to \cite{feireisl2019low}.

\subsection{The non-isentropic case}

We now consider non-isentropic fluids and take into account the transport of entropy. For the sake of simplicity, we study the Euler system for inhomogeneous fluids (we do not consider viscosity):
\begin{equation*}
\left\{
\begin{aligned}
    &\dt \rho + \div (\rho u) = 0,\\
    &\dt (\rho u) + \div (\rho u \otimes u) + \nabla p(\rho,S) =0,\\
    &\dt S + u \cdot \nabla S = 0,
\end{aligned}
\right.
\end{equation*}
where $S$ is the entropy and the pressure $p$ now depends on the density $\rho$ and the entropy. By convention, we write $\rho = R(p,S)$. We may for instance consider the pressure state law $p(\rho,S)=\rho^\gamma e^S$, $\gamma>1$. Let us note that the momentum equation can be rewritten as
\[
\rho \Big( \dt u+(u\cdot \nabla )u \Big) + \nabla p(\rho,S) = 0.
\]

By making the same change of variables than in the isentropic case \eqref{changevar} and scaling the entropy as $S(t,x) = S^\eps(\eps^2t,\eps x)$, a factor $1/\eps^2$ appears in front of the gradient of the pressure. Then, we follow the presentation of \cite{metivier2001limite} and make an adequate change of variable defining $q^\eps$ as $p^\eps(t,x) = \bar{p} e^{\eps q^\eps (\eps^2t,\eps x)}$ to obtain
\begin{equation*}
\left\{
\begin{aligned}
    &a^\eps (\partial_{t^\eps} q^\eps + u^\eps \cdot \nabla q^\eps) + \frac{1}{\eps} \div u^\eps =0,\\
   & r^\eps (\partial_{t^\eps} u^\eps + u^\eps \cdot \nabla u^\eps) + \frac{1}{\eps} \nabla q^\eps = 0,\\
   &\partial_{t^\eps} S^\eps + u^\eps \cdot \nabla S^\eps = 0,
\end{aligned}
\right.
\end{equation*}
where $a^\eps$ and $r^\eps$ are functions of $\rho^\eps$ and $S^\eps$ (see \cite{metivier2001limite} for more details). The initial data are $\rho^\eps_{\lvert t=0 }=\rho^\eps_0$, $u^\eps_{\lvert t=0 }=u^\eps_0$ and $S^\eps_{\lvert t=0 }=S^\eps_0$. We refer the interested reader to \cite{metivier2003averaging} for a proof of the existence of solutions of this system on a fixed time interval with uniform bounds in $\eps$.

Formally if we take the limit $\eps \to 0$, we obtain directly
\begin{align*}
    \div u^\eps \to 0,\\
    \nabla q^\eps \to 0.
\end{align*}

Therefore the pressures $q^0$ and thus $p^1$ are constant, where $p^1$ is the second term in the expansion $p^\eps=p^0+\eps p^1+\dots$. The limit of the pressure state law gives the limit density $\bar{\rho}$ as a function of the limit entropy $S$: $\bar{\rho} = R(\bar{p}, \bar{S})$. Setting $\nabla \pi^0$ as the limit of the pressure gradient, we obtain formally in the limit the Euler system for inhomogeneous incompressible fluids
\begin{equation*}
\left\{
    \begin{aligned}
        &r^0 (\dt u^0 + u^0 \cdot \nabla u^0) + \nabla \pi^0 = 0,\\
        &\dt S^0 + u^0 \cdot \nabla S^0 =0,\\
        &\div u^0 = 0. 
    \end{aligned}
\right.
\end{equation*}

The incompressible limit of non-isentropic Euler equations has been studied in the periodic case in \cite{metivier2003averaging} and in the whole space in \cite{metivier2001incompressible}. We also refer to \cite{schochet1986compressible} in the case of a bounded domain with well-prepared data and to \cite{alazard2005incompressible} in the general case. \cite{bresch2002low} also presents a multidimensional formal derivation of the low-Mach-number limit in an ill-prepared periodic case. For more details on the non-isentropic case, we refer the interested reader to the part of the non-isentropic equations presented in \cite{gallagher2005resultats} or to the presentation \cite{metivier2001limite}. For studies of the low-Mach-number limit in different types of compressible systems, we refer to \cite{lions1996mathematical} (in particular for the Navier-Stokes-Fourier system).

\section{Low-Mach number limit for two-phase flows with one single velocity}\label{sec2P1U}

We first investigate several one-velocity two-phase flow models, corresponding to different closure choices (in particular pressure equilibrium closure and volume fraction equation). Despite these differences, the presence of a single velocity field implies a unique characteristic flow velocity, and therefore a single Mach number. This setting allows us to introduce the low-Mach-number scaling in a unified manner and to analyze the resulting limit systems for the various closures considered. This assumption of one single velocity is physically justified in regimes where interfacial momentum exchange is sufficiently strong to enforce rapid velocity relaxation between phases. For small but non-zero relative velocity, we refer to the drift-flux system derived by \cite{ambroso2008drift,gastaldo2010scheme,gastaldo2011discretization}.

In this section and the next one, we focus on domains that are either periodic or bounded with the Dirichlet boundary conditions.

\subsection{Two-phase isentropic flows}\label{subE}

This subsection briefly reviews two-phase flows models (for more details, see the book \cite{ishii2010thermo} or the recent review in \cite{bresch2024mathematical} for instance).

First, one may mention as in \cite{gavrilyuk2020uncertainty} that there is no clear consensus on what constitutes a good two-phase flows model. Nevertheless, the majority of the models rely on the Baer–Nunziato equations, which describe compressible two-phase flows. One main difficulty is to understand and model both small scale and large scale effects and in particular the effects of small scale on larger scales. At a mesoscopic scale level, one may consider a two-phase flow as a flow composed of two one-phase areas separated by interfaces. Then, assuming some homogeneity at this mesoscale level and using tools as homogenization techniques, we can derive an averaged two-phase system that describes the flow at the macroscopic level. The interested reader may consult \cite{bresch2015note}, \cite{bresch2022mathematical} or \cite{bresch2023mathematical} to understand how averaged two-phase flow models may be obtained with no closure assumptions using homogenization techniques. Note that for the moment, the mathematical justification needs some restrictive assumptions like using the same velocity for both phases and only one dimension. Relaxing these assumptions is still a challenge. In the rest of this section, we will also consider only one velocity for both phases.

From a mathematical viewpoint, compressible two-fluid models have attracted increasing interest within the compressible PDE community; see, for instance, \cite{bresch2018multi}. These works clarify the structural properties of the systems and initiate their systematic analytical study.

Significant progress has been achieved in the viscous framework with weak regularity data. Finite-energy weak solutions were established in \cite{bresch2019finite} and \cite{vasseur2019global}, identifying suitable dissipation mechanisms. Further developments include weak and dissipative solution theories \cite{novotny2020weak,novotny2020weak1,jin2021existence}, weak–strong uniqueness results \cite{jin2019weak,li2021remarks,li2026weak}, and analysis of the inviscid case \cite{li2021weak}. Maximal regularity \cite{piasecki2022maximal}, large-time behavior \cite{li2020large}, and non-existence of classical solutions in the presence of vacuum \cite{li2025non} have also been investigated.

These contributions outline a coherent and evolving analytical theory for compressible two-phase flows.

\subsubsection{Two-phase systems with pressure algebraic closure.}

We now consider a fluid of two phases, $+$ and $-$, in a periodic or bounded domain $\Omega$ with the Dirichlet boundary conditions. We denote $\rho_i$, $p_i$, $\alpha_i$ the density, pressure and volume fraction of phase $i$. The single velocity is denoted by $u$ and any body force acting upon the fluid by $f$. Let $\mu_i$, $\lambda_i$ be the viscous coefficients of phase $i$ and $D(u)$ is the viscous tensor defined as
\[ D(u) = \frac{\nabla u + (\nabla u)^t}{2}.
\]

Different authors, such as \cite{novotny2020weak}, have studied a two-phase system with an algebraic closure for the pressure and only one velocity field taking the following form:
\begin{equation}
\left\{
    \begin{aligned}\label{algebraic1velocity}
        &\ap + \am = 1,\\
        &\dt (\ap \rp) + \div(\ap \rp u) = 0,\\
        &\dt (\am \rmi) + \div(\am \rmi u) = 0,\\
        &\dt ((\ap \rp + \am \rmi) u) + \div ((\ap \rp + \am \rmi) u \otimes u)\\
        &\quad \quad + \nabla (\ap p_+(\rp) + \am p_-(\rmi)) - 2 \div \Big( (\ap \mu_+(\rp) + \am \mu_-(\rmi)) D(u) \Big)\\
        &\quad \quad - \nabla \Big( (\ap \lambda_+ + \am \lambda_-) \div u \Big) = (\ap \rp + \am \rmi) f,\\
        &p_+(\rp) = p_-(\rmi).
    \end{aligned}
    \right.
\end{equation}

\begin{remark}
Note that the pressure term in the momentum equation can also be written as $\nabla (p_+(\rp))$ or $\ap \nabla p_+(\rp) + \am \nabla p_-(\rmi)$.
\end{remark} 

\begin{remark}
Note the link between $\div D(u)$ and $\Delta u$:
$$2 \mu \div D(u) = \mu \Delta u + \mu \nabla \div u. $$
Therefore, taking $\ap=1$, $\mu_+(\rp)=\mu$, $\lambda_+=\lambda$, we recover System \eqref{Sys1} for one-phase flow.\newline 
\end{remark}

\noindent \textbf{Formal energy equality.}

The energy equality associated with System \eqref{algebraic1velocity} reads for all $t \in[0, T]$:
\begin{multline}\label{energy_estimate_eps}
        \int_{\Omega} \frac{1}{2}\left(\alpha_+\rho_{+}+\alpha_-\rho_-\right)|u|^2(t,\cdot) +\int_{\Omega} \frac{\alpha_{+} (\rho_{+})^{\gamma_{+}}}{(\gamma_{+}-1)}(t,\cdot) +\int_{\Omega} \frac{\alpha_-(\rho_{-})^{\gamma_{-}}}{(\gamma_--1)}(t,\cdot) \\
        + \int_0^t\!\!\! \int_\Omega \Big( (\ap \mu_+(\rp) + \am \mu_-(\rmi)) \lvert D(u) \rvert^2 + (\ap \lambda_+ + \am \lambda_-) \lvert \div u \rvert^2 \Big) \\
        \quad \quad = \int_{\Omega} \frac{1}{2}\left(\alpha_{+} \rho_{+}+\alpha_-\rho_-\right)|u|^2(0,\cdot) +\int_{\Omega} \frac{\alpha_{+} (\rho_{+})^{\gamma_{+}}}{(\gamma_{+}-1)}(0,\cdot)+\int_{\Omega} \frac{\alpha_-(\rho_-)^{\gamma_-}}{(\gamma_{-}-1)}(0,\cdot) \\
        \quad \quad + \int_0^t\!\!\! \int_\Omega (\ap \rp + \am \rmi) f \cdot u.
\end{multline}

We provide for the reader's convenience the formal proof of this equality, when we assume that all variables are sufficiently regular. We start by taking the scalar product with $u$ of the momentum equation. Let's consider the first term:
\begin{align*}
    u \dt (\ap \rp u) &= \dt (\ap \rp |u|^2) - \ap \rp u \dt u,\\
     &= \dt (\ap \rp |u|^2) - u \Big( \dt (\ap \rp u) - u \dt (\ap \rp)\Big).
\end{align*}

Therefore we get:
\begin{align}\label{dtu}
    u \dt (\ap \rp u) = \frac{1}{2} \dt (\ap \rp |u|^2) + \frac{1}{2} |u|^2 \dt (\ap \rp).
\end{align}

In a similar way, we compute for the divergence term:
\begin{align}\label{divuu}
    u \div (\ap \rp u \otimes u) = \frac{1}{2} \div (\ap \rp |u|^2 u) + \frac{1}{2} |u|^2 \div (\ap \rp u).
\end{align}

Now adding \eqref{dtu} and \eqref{divuu} and using the mass equation for phase $+$ we obtain:
\begin{align}
    u \Big(\dt (\ap \rp u) + \div (\ap \rp u \otimes u) \Big) = \frac{1}{2} \dt (\ap \rp |u|^2) + \frac{1}{2} \div (\ap \rp |u|^2 u).\label{momu1}
\end{align}

Similarly we have for phase $-$:
\begin{align}
    u \Big(\dt (\am \rmi u) + \div (\am \rmi u \otimes u) \Big) = \frac{1}{2} \dt (\am \rmi |u|^2) + \frac{1}{2} \div (\am \rmi |u|^2 u).\label{momu2}
\end{align}

For the terms in $D(u)$ and $\div u$, we integrate by parts and under the assumption that the domain is either periodic or bounded with the Dirichlet boundary conditions, the boundary integrals vanish, such that we get:
\begin{align}
    - \int_\Omega \div \Big( (\ap \mu_+(\rp) + \am \mu_-(\rmi)) D(u) \Big) \cdot u - \int_\Omega \nabla \Big( (\ap \lambda_+ + \am \lambda_-) \div u \Big) \cdot u \nonumber \\
    = \int_\Omega \Big( (\ap \mu_+(\rp) + \am \mu_-(\rmi)) \lvert D(u) \rvert^2 + (\ap \lambda_+ + \am \lambda_-) \lvert \div u \rvert^2 \Big)\label{Duu}
\end{align}

Now we add \eqref{momu1}, \eqref{momu2} and \eqref{Duu} and integrate on the domain $\Omega$ the scalar product with $u$ of the momentum equation to obtain:
\begin{multline*}
    \frac{1}{2} \frac{\dd}{\dd t} \int_\Omega \left(\alpha_+\rho_{+}+\alpha_-\rho_-\right)|u|^2 +\\
    \int_\Omega \Big( (\ap \mu_+(\rp) + \am \mu_-(\rmi)) \lvert D(u) \rvert^2 + (\ap \lambda_+ + \am \lambda_-) \lvert \div u \rvert^2 \Big)\\
    + \int_\Omega \ap \nabla p_+(\rp) \cdot u + \int \am \nabla p_-(\rmi) \cdot u = \int_\Omega (\ap \rp + \am \rmi) f \cdot u.
\end{multline*}

The main novelty compared to the one-phase compressible flows is that there are two terms involving the pressures but also the volume fractions of each phase. We detail here the computations. Let's focus on the $+$ phase and set in this computation $\rho=\ap \rp$. Under the assumption that the domain $\Omega$ is either periodic or bounded with the Dirichlet boundary conditions, we obtain
\begin{align}
    \int_\Omega \ap \nabla p_+(\rp) \cdot u &= \int_\Omega \frac{\rho}{\rp} \nabla \rp^{\gp} \cdot u \nonumber \\
    &= \int_\Omega \frac{\gp}{\gp-1} \nabla \rp^{\gp-1} \cdot \rho u\nonumber\\
    &= \int_\Omega - \frac{\gp}{\gp-1} \rp^{\gp-1} \div(\rho u)\nonumber\\
    &= \int_\Omega \frac{\gp}{\gp-1} \rp^{\gp-1} \dt \rho\nonumber\\
    &= \int_\Omega \frac{\gp}{\gp-1} \rp^{\gp} \dt \ap + \int_\Omega \frac{\gp}{\gp-1} \ap \rp^{\gp-1} \dt \rp\nonumber\\
    &= \int_\Omega \frac{\gp}{\gp-1} \rp^{\gp} \dt \ap + \frac{\dd}{\dd t} \int_\Omega \frac{1}{\gp-1} (\ap \rp^{\gp}) - \int_\Omega \frac{1}{\gp-1} \rp^{\gp} \dt \ap\nonumber\\
    &= \frac{\dd}{\dd t} \int_\Omega \frac{1}{\gp-1} (\ap \rp^{\gp}) + \int_\Omega \rp^{\gp} \dt \ap.\label{pplus}\\ \nonumber
\end{align}

Similarly, we obtain for phase $-$
\begin{align}
    \int_\Omega \am \nabla p_-(\rmi) \cdot u = \frac{\dd}{\dd t} \int_\Omega \frac{1}{\gm-1} (\am \rmi^{\gm}) + \int_\Omega \rmi^{\gm} \dt \am.\label{pmoins}\\ \nonumber
\end{align}

We add \eqref{pplus} and \eqref{pmoins} using $\ap + \am =1$ and $\rp^{\gp} = \rmi^{\gm}$ we obtain the pressure terms in the energy inequality
\begin{align*}
    \frac{\dd}{\dd t} \int_\Omega \frac{\ap \rp^{\gp}}{\gp-1} + \frac{\am \rmi^{\gm}}{\gm-1}.
\end{align*}

There are few works concerning the well-posedness of such two-phase systems. We refer in particular to \cite{novotny2020weak} for the global existence of weak solution for averaged system with algebraic closure (when $\mu_+=\mu_-=\mu$ and $\lambda_+=\lambda_-=\lambda$).

\begin{remark} The energy \eqref{energie1} for one-phase systems is a particular case of this one, taking one volume fraction equal to 1 and the other to 0.\\
\end{remark}

\subsubsection{Two-phase systems with pressure PDE closure.}

Following \cite{bresch2024mathematical} we can consider a two-phase system with a PDE closure for the pressure and only one velocity field of the form:
\begin{equation}
    \left\{
        \begin{aligned}\label{PDE1velocity}
            &\ap + \am = 1,\\
            &\dt (\ap \rp) + \div(\ap \rp u) = 0,\\
            &\dt (\am \rmi) + \div(\ap \rmi u) = 0,\\
            &\dt ((\ap \rp + \am \rmi) u) + \div ((\ap \rp + \am \rmi) u \otimes u)\\
            &\quad \quad + \nabla \Big(\ap p_+(\rp) + \am p_-(\rmi)\Big)- 2 \div \Big( (\ap \mu_+(\rp) + \am \mu_-(\rmi)) D(u) \Big)\\
            &\quad \quad  - \nabla \Big( (\ap \lambda_+ + \am \lambda_-) \div u \Big) = (\ap \rp + \am \rmi) f,\\
            &\dt \ap + u \cdot \nabla \ap = \frac{1}{\eta} \ap \am (p_+(\rp)-p_-(\rmi)),
        \end{aligned}
    \right .
\end{equation}
with $\eta$ a relaxation coefficient.

\begin{remark}
    The only case in dimension $d>1$ for which the relaxation parameter can be explicitly characterized is the case of constant and equal viscosities, namely $\mu_+ = \mu_- = \mu$ and $\lambda_+ = \lambda_- = \lambda$, where $\mu$ and $\lambda$ are constants. In that case, the relaxation coefficient is given by $\eta = 2\mu + \lambda$.
    
    When $\mu_+$ and $\mu_-$ depend on $\rho$, the parameter $\eta$ should be viewed as an arbitrary relaxation coefficient. To the best of our knowledge, no rigorous justification is currently available, and this issue remains an open problem.
\end{remark}

\begin{remark}
    Note that when this relaxation coefficient tends to 0, we recover the case $p_+(\rp) = p_-(\rmi)$. On this subject, we refer to \cite{burtea2023pressure} and \cite{burtea2026relaxation} for more details.
\end{remark}

\noindent \textbf{Formal equality energy.}

When we assume that all variables are smooth enough, the associated energy equality reads for all $t \in[0, T]$:
\begin{multline}\label{energy_estimate_alpha}
        \int_{\Omega} \frac{1}{2}\left(\alpha_+ \rho_{+}+\alpha_-\rho_-\right)|u|^2(t,\cdot) +\int_{\Omega} \frac{\alpha_{+} (\rho_{+})^{\gamma_{+}}}{(\gamma_{+}-1)}(t,\cdot) +\int_{\Omega} \frac{\alpha_-(\rho_{-})^{\gamma_{-}}}{(\gamma_--1)}(t,\cdot) \\
        + \int_0^t\!\!\! \int_\Omega \Big( (\ap \mu_+(\rp) + \am \mu_-(\rmi)) \lvert D(u) \rvert^2 + (\ap \lambda_+ + \am \lambda_-) \lvert \div u \rvert^2 \Big) \\
        + \mbox{ \boldmath $\displaystyle \int_0^t\!\!\! \int_{\Omega} \frac{1}{\eta} \ap \am \Big(p_+(\rp)-p_-(\rmi)\Big)^2 $}\\
        = \int_{\Omega} \frac{1}{2}\left(\alpha_{+} \rho_{+}+\alpha_-\rho_-\right)|u|^2(0,\cdot) +\int_{\Omega} \frac{\alpha_{+} (\rho_{+})^{\gamma_{+}}}{(\gamma_{+}-1)}(0,\cdot) +\int_{\Omega} \frac{\alpha_-(\rho_-)^{\gamma_-}}{(\gamma_{-}-1)}(0,\cdot)\\
        +\int_0^t\!\!\! \int_\Omega (\ap \rp + \am \rmi) f \cdot u.
\end{multline}

Note the additional term in bold in Equation \eqref{energy_estimate_alpha} compared to the case with the algebraic closure for the pressure (Equation \eqref{energy_estimate_eps}). We detail here the computations for the terms involving the pressures. Let us focus on the $+$ phase. When we take the scalar product with $u$ of the momentum equation, we get the term $\nabla (\ap p_+(\rp)) \cdot u$ that we cannot simplify in the same way that we did when we had the algebraic closure. Here, we have:
\begin{align*}
    \int_\Omega \nabla (\ap p_+(\rp)) \cdot u = \int_\Omega p_+(\rp) u \cdot \nabla \ap + \int_\Omega \ap \nabla p_+(\rp) \cdot u
\end{align*}

For the first term we can use the PDE for the volume fraction and for the second term we can use \eqref{pplus}:
\begin{align*}
    \int_\Omega \nabla (\ap p_+(\rp)) \cdot u =& - \int_\Omega p_+(\rp) \dt \ap + \int_\Omega \frac{1}{\eta} \ap \am p_+(\rp) \Big(p_+(\rp)-p_-(\rmi)\Big)\\
    &+ \frac{\dd}{\dd t} \int_\Omega \frac{1}{\gp-1} (\ap \rp^{\gp}) + \int_\Omega p_+(\rp) \dt \ap\\
    =& \int_\Omega \frac{1}{\eta} \ap \am p_+(\rp) \Big(p_+(\rp)-p_-(\rmi)\Big)+ \frac{\dd}{\dd t} \int_\Omega \frac{1}{\gp-1} (\ap \rp^{\gp}).\\
\end{align*}

As $\partial_t \am + u\cdot \nabla \am = -(\partial_t \ap + u\cdot \nabla \ap)$, we obtain similarly for phase $-$
\begin{align*}
    \int_\Omega \nabla (\am p_-(\rmi)) \cdot u =& - \int_\Omega \frac{1}{\eta} \ap \am p_-(\rmi) \Big(p_+(\rp)-p_-(\rmi)\Big) + \frac{\dd}{\dd t} \int_\Omega \frac{1}{\gm-1} (\am \rmi^{\gm}).\\
\end{align*}

Therefore, adding the two terms, we obtain
\begin{multline*}
    \int_\Omega \nabla(\ap p_+(\rp) + \am p_-(\rmi)) \cdot u = \int_\Omega \frac{1}{\eta} \ap \am \Big(p_+(\rp) - p_-(\rmi) \Big)^2 \\
    + \frac{\dd}{\dd t} \int_\Omega  \frac{\ap \rp^{\gp}}{\gp-1}+\frac{\am \rmi^{\gm}}{\gm-1}.
\end{multline*}

Recently, \cite{novotny2020weak1} has studied a system similar to System \eqref{PDE1velocity}, considering constant viscosities and a transport equation for the volume fraction (no relaxation term in the PDE closure) because he considers a diffuse interface. Note that this makes his system not exactly an averaged mixture system. For this system, he proves the global existence of weak solutions.

In the rest of this section, for the sake of mathematical simplicity, we make some additional assumptions: we set $f=0$ and use the same viscosity coefficients $\mu_+(\rp)=\mu_-(\rmi)=\mu$ and $ \lambda_+=\lambda_-=\lambda$, with $\mu$ and $\lambda$ constant. We emphasize that the more general case of phase-dependent viscosity coefficients $\mu_+(\rho_+)$, $\mu_-(\rho_-)$ and distinct bulk viscosities $\lambda_+$, $\lambda_-$ remains an important and challenging problem.

\subsubsection{Formal derivation of the low-Mach-number limit for two-phase systems with pressure algebraic closure}\label{formalp}

In this subsection, we consider System \eqref{algebraic1velocity} and we do a dimensional analysis as explained in the first section of this paper to make appear the Mach number in the equations (change of variable \eqref{changevar}).
We consider the following a-dimensional two-phase flow system.
\begin{align}
    \left\{
        \begin{array}{l}\label{system_eps}
            \partial_t (\alpha_+^{\epsilon} \rho_+^{\epsilon}) + \div(\alpha_+^{\epsilon}\rho_+^{\epsilon}u^{\epsilon}) = 0,\\
            \partial_t (\alpha_-^{\epsilon} \rho_-^{\epsilon}) + \div(\alpha_-^{\epsilon}\rho_-^{\epsilon}u^{\epsilon}) = 0,\\
            \partial_t ((\alpha_+^{\epsilon} \rho_+^{\epsilon} + \alpha_-^{\epsilon}\rho_-^{\epsilon})u^{\epsilon}) + \div((\alpha_+^{\epsilon}\rho_+^{\epsilon} + \alpha_-^{\epsilon}\rho_-^{\epsilon})u^{\epsilon}\otimes u^{\epsilon})\\
            \quad \quad + \frac{1}{\epsilon^2} \nabla (\alpha_+^{\epsilon}p_+(\rho_+^{\epsilon})+\alpha_-^{\epsilon}p_-(\rho_-^{\epsilon})) = \mu \Delta u^{\epsilon} + (\mu + \lambda)\nabla \div u^{\epsilon},\\
            p_i(\rho_i^{\epsilon}) = (\rho_i^{\epsilon})^{\gamma_i}, i=+,-\\
            p_+(\rho_+^{\epsilon}) = p_-(\rho_-^{\epsilon}),\\
            \alpha_+^{\epsilon} + \alpha_-^{\epsilon} = 1.
        \end{array}
    \right .
\end{align}
For the sake of readability, we have dropped, here and in the sequel, the subscript $\eps$ for the time and space variables, writing $(t,x)$ instead to $(t^\eps,x^\eps)$.

In this part, we will give a formal proof of the following proposition.

\begin{proposition}
We consider a smooth enough solution of system \eqref{system_eps} associated with well-prepared initial data in the sense of Definition \eqref{defwpd}, in a domain that is either periodic or bounded with the Dirichlet boundary conditions, meaning $u = 0$ on $\partial \Omega$. Taking the limit as $\eps$ goes to zero in the system \eqref{system_eps}, we get formally the following limit system:
\begin{equation}\label{limitsys}
    \begin{cases}
        \rp^0 = C^{\frac{1}{\gp}},\\
        \rmi^0 = C^{\frac{1}{\gm}},\\
        \div u^0 = 0,\\
        \dt(\left(\rp^0 \ap^0 + \rmi^0 \am^0 \right) u^0) + \div( \left(\rp^0 \ap^0 + \rmi^0 \am^0 \right) u^0\otimes u^0) + \nabla \pi^0 = \mu \Delta u^0,\\
        \dt \ap^0 + u^0 \cdot \nabla \ap^0 =0,\\
        \ap^0 + \am^0 = 1,
    \end{cases}
\end{equation}
with $\nabla \pi^0 = p_+'(\rho_+^0) \nabla \rho_+^1$, being the ``limit'' of $\frac{(\rho^\eps)^{\gamma_+}-(\rho^0)^{\gamma_+}}{\eps^2}=\frac{p_+(\rho^\eps)-p_+(\rho^0)}{\eps^2}$.
\end{proposition}

\begin{remark}
Note that the transport equation on $\alpha_+^0$ also implies a transport equation on $\rho^0=\rp^0 \ap^0 + \rmi^0 \am^0$ as $\rp^0$ and $\rmi^0$ are constant:
\begin{align*}
    \dt \rho^0 + u^0 \cdot \nabla \rho^0 = 0.
\end{align*}
and the momentum equation can also be rewritten on the form
\[
\rho^0 \left( \dt u^0 +  u^0\cdot \nabla  u^0 \right) + \nabla \pi^0 = \mu \Delta u^0.
\]
\end{remark}

\begin{remark}
    The rigorous justification of this incompressible limit can be found in \cite{lebot2026low}.
\end{remark}

\begin{remark}
For the sake of simplicity, in the rest of the paper we consider that $\ap$ is not zero on a set of non-zero measure.
\end{remark}

\begin{proof} [Formal proof.]

To obtain the asymptotics, we expand the solutions with respect to $\eps$ as follows
\begin{align}
    \begin{array}{cc}\label{expansion}
        \rho_i^\eps =& \rho_i^0 + \eps^2 \rho_i^1 + \eps^4 \rho_i^2,\\
        \alpha_i^\eps =& \alpha_i^0 + \eps^2 \alpha_i^1 + \eps^4 \alpha_i^2,\\
        u^\eps =& u^0 + \eps^2 u^1 + \eps^4 u^2,\\
    \end{array}
\end{align}
with $\eps$ sufficiently small. Let us notice that the only small parameter in System \eqref{system_eps} is $\eps^2$, so it is natural to introduce an expansion in terms of power of $\eps^2$. Alternatively, we could have set $\tilde \epsilon:=\eps^2$ and used an expansion in terms of power of $\tilde \eps$.

This implies the following pressure expansion with respect to $\eps$
\begin{align}
    \begin{array}{cc}\label{expansionp}
        p_+(\rho_+^\eps) &= p_+(\rho_+^0) + \eps^2 p_+'(\rho_+^0) (\rho_+^1 + \eps^2 \rho_+^2+...),\\
        p_-(\rho_-^\eps) &= p_-(\rho_-^0) + \eps^2 p_-'(\rho_-^0)(\rho_-^1 + \eps^2 \rho_-^2+...).
    \end{array}
\end{align}

At order $1$ and order $\eps^2$, the algebraic constraint $\ap^\eps + \am^\eps =1$ reads respectively
\begin{align}
    \alpha_+^0 + \alpha_-^0 = 1, \label{alpha0}\\
    \ap^1 + \am^1 = 0.\label{alpha1}
\end{align}

The equality of pressures gives, respectively at order $1$ and $\eps^2$:
\begin{align}
    p_+(\rho_+^0) &= p_-(\rho_-^0),\label{p0} \\
    p_+'(\rho_+^0) \rho_+^1 &= p_-'(\rho_-^0)\rho_-^1. \label{p1}
\end{align}

At order $\eps^{-2}$, the momentum equation reads
\begin{align*}
    \nabla (\alpha_+^{0}p_+(\rho_+^0)+\alpha_-^0p_-(\rho_-^0)) =0.\\ \nonumber
\end{align*}

By virtue of \eqref{alpha0} and \eqref{p0}, this last equation yields
\begin{align*}
    \left\{
        \begin{array}{l}
              \nabla p_+(\rho_+^0) = 0, \\
              \nabla p_-(\rho_-^0) = 0.
        \end{array}
    \right.
\end{align*}

Let's call $C$ a constant in space, from the equality of pressures \eqref{p0} we deduce:
\begin{align}
    \rho_+^0 = C(t)^{1/\gamma_+},\label{rho+0}\\
    \rho_-^0 = C(t)^{1/\gamma_-}.\label{rho-0}
\end{align}

The mass equations expansion yields at order $1$:
\begin{align}
    \dt (\alpha_+^0 \rho_+^0) + \div (\alpha_+^0 \rho_+^0 u^0) = 0,\label{mass+0}\\
    \dt (\alpha_-^0 \rho_-^0) + \div (\alpha_-^0 \rho_-^0 u^0) = 0.\label{mass-0}
\end{align}

Plugging \eqref{rho+0} and \eqref{rho-0} into \eqref{mass+0} and \eqref{mass-0} respectively, we obtain
\begin{align*}
    \alpha_+^0 \frac{\dt C^{1/\gamma_+}}{C^{1/\gamma_+}} + \dt \alpha_+^0 + \div(\alpha_+^0 u^0) = 0,\\
    \alpha_-^0 \frac{\dt C^{1/\gamma_-}}{C^{1/\gamma_-}} + \dt \alpha_-^0 + \div(\alpha_-^0 u^0) = 0.
\end{align*}

By summing these two last equations, using $\alpha_+^0 + \alpha_-^0 =1$, we obtain
\begin{align}
    \alpha_+^0 \dt \ln (C^{1/\gamma_+}) + \alpha_-^0 \dt \ln(C^{1/\gamma_-}) + \div u^0 = 0.\label{fordivu}
\end{align}

We now integrate in space over $\Omega$, that we recall is considered periodic or bounded with the Dirichlet boundary conditions, we get:
\begin{align*}
    &\frac{1}{\lvert \Omega \rvert} \left[ \int_\Omega \Big(  \alpha_+^0  \Big) \dt \ln (C^{1/\gamma_+}) + \int_\Omega \Big( \alpha_-^0 \Big) \dt \ln(C^{1/\gamma_-}) \right] = 0\\
    &\Rightarrow \left( \frac{1}{\gamma_+} \int_\Omega \alpha_+^0 + \frac{1}{\gamma_-} \int_\Omega \alpha_-^0 \right) \frac{\dt C}{C} = 0.
\end{align*}

Therefore, if $\displaystyle \frac{1}{\gamma_+} \int_\Omega \alpha_+^0 + \frac{1}{\gamma_-} \int_\Omega \alpha_-^0 \neq 0$, which is always true by positivity of the volume fractions, $\dt C = 0$ and we get, with $C\rvert_{t=0} = C_0$
\begin{align*}
    \rho_+^0 = C_0^{1/\gamma_+},\\
    \rho_-^0 = C_0^{1/\gamma_-}.
\end{align*}

Using \eqref{fordivu}, we obtain the divergence free criterion:
\begin{align*}
    \div u^0 =0.
\end{align*}

From the order 1 of the mass equations \eqref{mass+0} or \eqref{mass-0}, we get
\begin{align}
    \dt \alpha_+^0 + u^0 \cdot \nabla \alpha_+^0 =0.\label{eqap}
\end{align}

Finally, using the order 1 of the momentum equation:
\begin{multline*}
    \dt ( (\rho_+^0 \alpha_+^0 + \alpha_-^0 \rho_-^0) u^0) + \div \Big( (\rho_+^0 \alpha_+^0 + \alpha_-^0 \rho_-^0) u^0 \otimes u^0  \Big)\\
    + \nabla \Big( \ap^0 p_+'(\rho_+^0) \rho_+^1 + \am^0 p_-'(\rho_-^0)\rho_-^1 + \ap^1 p_+(\rp^0) + \am^1 p_-(\rmi^0) \Big) = \mu \Delta u^0.
\end{multline*}

Using \eqref{alpha0}, \eqref{alpha1}, \eqref{p0} and \eqref{p1} and the fact that $\rp^0$ and $\rmi^0$ are constant, we get:
\begin{align*}
     \dt ( (\rho_+^0 \alpha_+^0 + \alpha_-^0 \rho_-^0) u^0) + \div \Big( (\rho_+^0 \alpha_+^0 + \alpha_-^0 \rho_-^0) u^0 \otimes u^0  \Big) + \nabla \Big(  p_+'(\rho_+^0) \rho_+^1 \Big) = \mu \Delta u^0.
\end{align*}

Setting $\pi^0 :=  p_+'(\rho_+^0) \rho_+^1$ and using \eqref{eqap}, we can also write:
\begin{align*}
    (\rho_+^0 \alpha_+^0 + \alpha_-^0 \rho_-^0) \left( \dt u^0 + u^0 \cdot \nabla u^0 \right) + \nabla \pi^0 - \mu \Delta u^0 =0.
\end{align*}

We have shown formally that the limit system obtained when $\eps \to 0$ is \eqref{limitsys}.\\
\end{proof}

\subsubsection{Formal derivation of the low-Mach-number limit for two-phase systems with pressure PDE closure}

In this part, we consider the following system with a PDE closure for the pressure:
\begin{equation}\label{sys_alpha_eps}
    \left\{
        \begin{aligned}
            &\partial_t (\alpha_+^{\epsilon} \rho_+^{\epsilon}) + \div(\alpha_+^{\epsilon}\rho_+^{\epsilon}u^{\epsilon}) = 0,\\
            &\partial_t (\alpha_-^{\epsilon} \rho_-^{\epsilon}) + \div(\alpha_-^{\epsilon}\rho_-^{\epsilon}u^{\epsilon}) = 0,\\
            &\partial_t ((\alpha_+^{\epsilon} \rho_+^{\epsilon} + \alpha_-^{\epsilon}\rho_-^{\epsilon})u^{\epsilon}) + \div((\alpha_+^{\epsilon}\rho_+^{\epsilon} + \alpha_-^{\epsilon}\rho_-^{\epsilon})u^{\epsilon}\otimes u^{\epsilon})\\
            & \quad \quad + \frac{1}{\epsilon^2} \nabla (\alpha_+^{\epsilon}p_+(\rho_+^{\epsilon})+\alpha_-^{\epsilon}p_-(\rho_-^{\epsilon})) = \mu \Delta u^{\epsilon} + (\mu + \lambda)\nabla \div u^{\epsilon},\\
            &p_i(\rho_i^{\epsilon}) = (\rho_i^{\epsilon})^{\gamma_i}, i=+,-,\\
            &\dt \ap^\eps + u^\eps \cdot \nabla \ap^\eps = \frac{1}{\eps^2 \eta} \ap^\eps \am^\eps \left(p_+(\rp^\eps)-p_-(\rmi^\eps)\right),\\
            &\alpha_+^{\epsilon} + \alpha_-^{\epsilon} = 1.
        \end{aligned}
    \right .
\end{equation}

We will give a formal proof of the following proposition.

\begin{proposition}
We consider a smooth enough solution of system \eqref{sys_alpha_eps} associated with well-prepared initial data in the sense of Definition \eqref{defwpd}, in a domain that is either periodic or bounded with the Dirichlet boundary conditions, meaning $u = 0$ on $\partial \Omega$. Formally taking the low-Mach-number limit of \eqref{sys_alpha_eps} then yields system \eqref{limitsys}.
\end{proposition}

\begin{proof} [Formal proof.]

To obtain the asymptotics for this system, we use the same expansions \eqref{expansion} and \eqref{expansionp} as for the system with the algebraic closure for the pressure.

At order 1, the algebraic constraint on $\ap^\eps$ also reads
\begin{align*}
    \alpha_+^0 + \alpha_-^0 = 1. 
\end{align*}

But now we obtain the equality of pressures from the PDE closure at order $\eps^{-2}$:
\begin{align*}
    \ap^0 \am^0 (p_+(\rp^0) - p_-(\rmi^0)) =0 \Rightarrow p_+(\rp^0) = p_-(\rmi^0).
\end{align*}

Note that we are not interested in the case where $\ap =0$ or $\am=0$ as it amounts to consider only one-phase flow.\\

The momentum equation reads at order $\eps^{-2}$
\begin{align*}
    \nabla (\alpha_+^{0}p_+(\rho_+^0)+\alpha_-^0p_-(\rho_-^0)) =0. \nonumber
\end{align*}

In a same way than in the case of the algebraic closure, we can deduce, using the equality of pressures
\begin{align}
    \rho_+^0 = C(t)^{1/\gamma_+},\label{rp0}\\
    \rho_-^0 = C(t)^{1/\gamma_-}.\label{rmi0}
\end{align}
with $C$ a constant in space.

The mass equations read at order $1$:
\begin{align*}
    \dt (\alpha_+^0 \rho_+^0) + \div (\alpha_+^0 \rho_+^0 u^0) = 0,\\
    \dt (\alpha_-^0 \rho_-^0) + \div (\alpha_-^0 \rho_-^0 u^0) = 0.
\end{align*}

Summing these two equations and using \eqref{rp0}-\eqref{rmi0}, we get:
\begin{align}
    \alpha_+^0 \dt \ln (C^{1/\gamma_+}) + \alpha_-^0 \dt \ln(C^{1/\gamma_-}) + \div u^0 = 0.\label{fordivu2}
\end{align}

Similarly than in the first case, we integrate in space and obtain that $C$ is a constant and
\begin{align*}
    \rho_+^0 = C_0^{1/\gamma_+},\\
    \rho_-^0 = C_0^{1/\gamma_-}.
\end{align*}

Then from \eqref{fordivu2} we get $\div u^0 =0$ and from the order 1 of the mass equations we get 
\begin{align}
    \dt \ap^0 + u^0 \cdot \nabla \ap^0 =0.\label{ap0}
\end{align}

It remains to derive the momentum equation. For that, we use the order $1$ of the PDE closure:
\begin{align*}
    \dt \ap^0 + u^0 \cdot \nabla \ap^0 = \frac{1}{\mu} \ap^0 \am^0 \Big( p_+'(\rp^0)\rp^1 - p_-'(\rmi^0)\rmi^1 \Big)
\end{align*}

From \eqref{ap0}, we deduce that $p_+'(\rp^0)\rp^1 = p_-'(\rmi^0)\rmi^1 := \pi^0$ and we can plug it in the order~$1$ of the momentum equation to get
\begin{align*}
    \dt(\left(\rp^0 \ap^0 + \rmi^0 \am^0 \right) u^0) + \div( \left(\rp^0 \ap^0 + \rmi^0 \am^0 \right) u^0\otimes u^0) + \nabla \pi^0 = \mu \Delta u^0.
\end{align*}

We have obtained formally the same limit system \eqref{limitsys} when $\eps$ goes to 0 than in the case of an algebraic closure.
\end{proof}

\subsection{Two-phase non-isentropic flows}

We consider the following simplified two-phase non-isentropic flow system, with a pressure algebraic closure and where $S_\pm$ is the entropy of phase $+/-$. This system reads:
\begin{equation*}
    \left\{
        \begin{aligned}
            &\partial_t (\alpha_+ \rho_+) + \div(\alpha_+\rho_+u) = 0,\\
            &\partial_t (\alpha_- \rho_-) + \div(\alpha_-\rho_-u) = 0,\\
            &\partial_t ((\alpha_+ \rho_+ + \alpha_-\rho_-)u) + \div((\alpha_+\rho_+ + \alpha_-\rho_-)u\otimes u) \\
            & \quad \quad + \nabla (\alpha_+p_+(\rho_+,S_+)+\alpha_- p_-(\rho_-,S_-)) = \mu \Delta u + (\mu + \lambda)\nabla \div u,\\
            &p_i(\rho_i,S_i) = (\rho_i)^{\gamma_i}e^{S_i}, i=+,-\\
            &p_+(\rho_+,S_+) = p_-(\rho_-,S_-),\\
            &\alpha_+ + \alpha_- = 1.,\\
            &\dt S_\pm + u \cdot \nabla S_\pm = 0.
        \end{aligned}
    \right.
\end{equation*}

Note that now the pressure of one phase depends on its density and its entropy.

After the a-dimensionalization, we obtain
\begin{equation}\label{noni1veleps}
    \begin{cases}
            \partial_t (\alpha_+^{\epsilon} \rho_+^{\epsilon}) + \div(\alpha_+^{\epsilon}\rho_+^{\epsilon}u^{\epsilon}) = 0,\\
            \partial_t (\alpha_-^{\epsilon} \rho_-^{\epsilon}) + \div(\alpha_-^{\epsilon}\rho_-^{\epsilon}u^{\epsilon}) = 0,\\
            \partial_t ((\alpha_+^{\epsilon} \rho_+^{\epsilon} + \alpha_-^{\epsilon}\rho_-^{\epsilon})u^{\epsilon}) + \div((\alpha_+^{\epsilon}\rho_+^{\epsilon} + \alpha_-^{\epsilon}\rho_-^{\epsilon})u^{\epsilon}\otimes u^{\epsilon})\\
            \quad \quad + \frac{1}{\epsilon^2} \nabla (\alpha_+^{\epsilon}p_+(\rho_+^{\epsilon},S_+^\eps)+\alpha_-^{\epsilon}p_-(\rho_-^{\epsilon},S_-^\eps)) = \mu \Delta u^{\epsilon} + (\mu + \lambda)\nabla \div u^{\epsilon},\\
            p_i(\rho_i^{\epsilon},S_i^\eps) = (\rho_i^{\epsilon})^{\gamma_i} e^{S_i^\eps}, i=+,-\\
            p_+(\rho_+^{\epsilon},S_+^\eps) = p_-(\rho_-^{\epsilon},S_-^\eps),\\
            \alpha_+^{\epsilon} + \alpha_-^{\epsilon} = 1,\\
            \dt S_\pm^\eps + u^\eps \cdot \nabla S_\pm^\eps = 0.
    \end{cases}
\end{equation}

In this part, we will give a formal proof of the following proposition.

\begin{proposition} 
We consider a smooth enough solution of system \eqref{noni1veleps} associated with well-prepared initial data in the sense of Definition \eqref{defwpd}, in a domain that is either periodic or bounded with the Dirichlet boundary conditions. Taking the limit as $\eps$ goes to zero in the system \eqref{noni1veleps}, we get formally the following limit system:
\begin{equation*}
    \begin{cases}
        \rp^0 = C^{\frac{1}{\gp}} e^{-S_+/\gp},\\
        \rmi^0 = C^{\frac{1}{\gm}} e^{-S_-/\gm},\\
        \div u^0 = 0,\\
         \dt \left( \left(\rp^0 \ap^0 + \rmi^0 \am^0 \right) u^0\right) + \div \Big( \left(\rp^0 \ap^0 + \rmi^0 \am^0 \right) u^0\otimes  u^0\Big) + \nabla \pi^0 = \mu \Delta u^0,\\
        \dt \ap^0 + u^0 \cdot \nabla \ap^0 =0,\\
        \ap^0 + \am^0 = 1,\\
        \dt S_\pm^0 + u^0 \cdot \nabla S_\pm^0=0.
    \end{cases}
\end{equation*}
with $\pi^0$ the limit pressure.
\end{proposition}

\begin{remark}
Note that because the pressure now depends on the density and the entropy, the density is not constant anymore in the limit. It follows a transport equation: 
\begin{align*}
    \dt \rho_\pm^0 + u^0 \cdot \nabla \rho_\pm^0 =0.
\end{align*}
\end{remark}

\begin{proof}[Formal proof]
We use the same expansions as before, except for the pressures that now depend on the density and the entropy. Its expansion becomes:
\begin{align}
    p_i(\rho_i,S_i) = p_i(\rho_i^0,S_i^0) + \eps^2 \Big( \partial_\rho p_i(\rho_i^0,S_i^0) \rho_i^1 + \partial_S p_i(\rho_i^0,S_i^0) S_i^1 + \cdots \Big).\label{expp}
\end{align}

At order 1, we have
\begin{align}
    \ap^0 + \am^0 = 1,\label{apm0}\\
    p_+ (\rp^0,S_+^0) = p_-(\rmi^0,S_-^0),\label{eqp}\\
    \dt S_\pm^0 + u_\pm^0 \cdot \nabla S_\pm^0 = 0.\label{entropy0}
\end{align}

Using the two first equalities, we get that at order $\eps^{-2}$, the momentum equations reads
\begin{align*}
    \nabla p_+ (\rp^0,S_+^0) = 0,\\
    \nabla p_-(\rmi^0,S_-^0)=0.
\end{align*}

Therefore, setting $C(t)$ a constant in space, from Equation \eqref{eqp} we have
\begin{align*}
     p_+ (\rp^0,S_+^0) = p_-(\rmi^0,S_-^0) = C(t).
\end{align*}

From the state law, we can deduce that
\begin{align*}
    \rho_\pm^0 = C(t)^{1/\gamma_\pm} e^{-S_\pm^0/\gamma_\pm}.
\end{align*}

Using the mass equations at order 1 and \eqref{entropy0}, we get
\begin{align*}
    \rho_\pm^0(\dt \alpha_\pm^0 + \div(\alpha_\pm^0 u^0)) = - \alpha_\pm^0 (\dt \rho_\pm^0 + u^0 \cdot \nabla \rho_\pm^0) = - \frac{\alpha_\pm^0 \rho_\pm^0}{\gamma_\pm} \frac{C'(t)}{C(t)}.
\end{align*}

Adding these two equations and using \eqref{apm0}, we have
\begin{align}
    \div (u^0) + \Big( \frac{\ap^0}{\gp} + \frac{\am^0}{\gm} \Big) \frac{C'(t)}{C(t)}=0.\label{divuI}
\end{align}

Now we integrate in space. Because we have assumed a periodic domain or a bounded domain with the Dirichlet boundary conditions, we get
\begin{align*}
    \frac{C'(t)}{C(t)} \int_\Omega \Big( \frac{\ap^0}{\gp} + \frac{\am^0}{\gm} \Big) d x   =0.
\end{align*}

Thus $C'(t) = 0$ and $C$ is a constant in time and space. Coming back to \eqref{divuI} we obtain
\begin{align*}
    \div u^0 = 0.
\end{align*}

Now from the order 1 of the mass equations, we get
\begin{align*}
    \dt \ap^0 + u^0 \cdot \nabla \ap^0 =0.
\end{align*}

Finally we can derive the limit momentum equation. For that, we use the expansion of the pressure \eqref{expp} and the order $\eps$ of the equality of pressures:
\begin{align*}
    \Big( \gamma_+ \frac{\rp^1}{\rp^0} + S_+^1 \Big) p_+(\rp^0,S_+^0) = \Big( \gm \frac{\rmi^1}{\rmi^0} + S_-^1 \Big) p_-(\rmi^0,S_-^0).
\end{align*}

We define
\begin{align*}
    \pi_i^0 = \Big( \gamma_i \frac{\rho_i^1}{\rho_i^0} + S_i^1 \Big) p_i(\rho_i^0,S_i^0).
\end{align*}

Hence $\pi_+^0 = \pi_-^0$, and the order 1 of the momentum equation reads:
\begin{align*}
    \left(\rp^0 \ap^0 + \rmi^0 \am^0 \right) \left( \dt u^0 + \div (u^0 \otimes u^0)\right) + \nabla \Big(\ap^0 \pi_+^0 + \am^0 \pi_-^0  \Big) = \mu \Delta u^0, \nonumber \\
    \left(\rp^0 \ap^0 + \rmi^0 \am^0 \right) \left( \dt u^0 + \div (u^0 \otimes u^0)\right) + \nabla \pi_+^0 = \mu \Delta u^0.
\end{align*}
\end{proof}

\section{Low-Mach number limit for two-phase two-velocities flows}\label{sec2P2U}

The second part of this work is devoted to two-velocity two-phase flow models. Although the phases may have distinct velocities, we assume that these velocities remain of the same order. This allows us to introduce a single Mach number and to consider a common low-Mach-number regime. This weakly compressible asymptotic setting corresponds to flows close to mechanical equilibrium, in which interfacial coupling ensures a rapid equilibration of pressure. As a result, acoustic effects do not contribute at leading order, and a single low-Mach asymptotic framework can be used to formally derive the limit system. The following analysis is formal and aims at identifying the expected structure of the low-Mach limit for a class of two-velocity two-phase models. Let us emphasize that allowing for distinct Mach numbers for each phase would lead to a significantly more intricate asymptotic structure, as different compressibility scales would interact. Such regimes require a separate analysis and fall beyond the scope of the present work; we refer to the existing literature for related investigations (\cite{benjelloun2021sound}).

We recall that we focus on domains that are either periodic or bounded with appropriate boundary conditions.

\subsection{The isentropic case}

In this part, we focus on the isentropic case with again a barotropic state law $p_\pm=p_\pm(\rho_\pm)$. Each phase $\pm$ is now also characterised by its velocity field $u_\pm$. Following \cite{bresch2024mathematical}, we consider the two-phase two-velocities system with a pressure algebraic closure:
\begin{align*}
    \left\{
        \begin{array}{l}
            \ap + \am = 1,\\
            \dt (\ap \rp) + \div(\ap \rp u_+) = 0,\\
            \dt (\am \rmi) + \div(\ap \rmi u_-) = 0,\\
            \dt (\ap \rp u_+) + \div (\ap \rp u_+ \otimes u_+) - 2 \div(\ap \mu_+(\rp) D(u_+))\\ 
            \quad \quad- \nabla ( \ap \lambda_+ (\rp) \div u_+ ) + \ap \nabla ( p_+(\rp)) + p_{\mathrm{int},+} \nabla \ap = \ap \rp f + F_{\mathrm{int},+},\\
            \dt (\am \rmi u_-) + \div (\am \rmi u_- \otimes u_-) - 2 \div(\am \mu_-(\rmi) D(u_-))\\
            \quad \quad - \nabla ( \am \lambda_- (\rmi) \div u_- ) + \am \nabla ( p_-(\rmi)) + p_{\mathrm{int},-} \nabla \am = \am \rmi f  + F_{\mathrm{int},-},\\
            p = p_+(\rp) = p_-(\rmi),\\
            (\ap \rp)\rvert_{t=0} = R_{+,0}, \quad (\ap \rp u_+)\rvert_{t=0} = m_{+,0},\\
            (\am \rmi)\rvert_{t=0} = R_{-,0}, \quad (\am \rmi u_-)\rvert_{t=0} = m_{-,0}.
        \end{array}
    \right.
\end{align*}
with $p_{\mathrm{int},i}$ and $F_{\mathrm{int},i}$ respectively the pressure and the drag term of phase $i=+,-$ at the interface. 

\begin{remark}
We refer the interested reader to \cite{bresch2018multi} for a discussion on the term $p_{\mathrm{int},i}$. Note that it is an important choice for the model because it ensures the hyperbolicity of the system. Here we will assume that $p_{\mathrm{int},+}=p_{\mathrm{int},-}=p_{\mathrm{int}}$.
\end{remark}

\begin{remark}
We can specify the drag terms as
$$F_{\mathrm{int},\pm} = \frac{1}{\eta} \lvert u_+ - u_- \rvert (u_\pm-u_\mp). $$
Note that if we take the limit $\eta \to 0$, we get $u_+=u_-$ and by adding the two momentum equations we recover System \eqref{algebraic1velocity}.
\end{remark}

There are few works on the well-posedness of such system. We note in particular the work of \cite{bresch2012global} that shows the global existence of weak solution in the one-dimensional case and under some assumptions on the viscosity coefficients and the pressure laws. In the multidimensional case we refer to \cite{bresch2010global} (that includes capillarity in the system). 

In the following, for the sake of simplicity, we do not consider any viscosity effects nor surface tension and we consider only the interfacial pressure. We do the adimensional analysis of the system using the same Mach number for both phases, with the same change of variable \eqref{changevar}. It leads to a factor of $1/\eps^2$ in front of the pressure terms in both momentum equations.
\begin{align}
    \left\{
        \begin{array}{l}\label{algebraic2velocityeps}
            \ap^\eps + \am^\eps = 1,\\
            \dt (\ap^\eps \rp^\eps) + \div(\ap^\eps \rp^\eps u_+^\eps) = 0,\\
            \dt (\am^\eps \rmi^\eps) + \div(\ap^\eps \rmi^\eps u_-^\eps) = 0,\\
            \dt (\ap^\eps \rp^\eps u_+^\eps) + \div (\ap^\eps \rp^\eps u_+^\eps \otimes u_+^\eps) + \frac{1}{\eps^2} \ap^\eps \nabla p_+(\rp^\eps) +p_{\mathrm{int}}^\eps \nabla \ap^\eps=0,\\
            \dt (\am^\eps \rmi^\eps u_-^\eps) + \div (\am^\eps \rmi^\eps u_-^\eps \otimes u_-^\eps) + \frac{1}{\eps^2} \am^\eps \nabla p_-(\rmi^\eps) + p_{\mathrm{int}}^\eps \nabla \am^\eps=0,\\
            p_+(\rp^\eps) = p_-(\rmi^\eps).
        \end{array}
    \right.
\end{align}

In the following, we are going to give a formal proof of the following proposition.

\begin{proposition} 
Consider a smooth enough solution of system \eqref{algebraic2velocityeps} associated with well-prepared initial data in the sense of Definition \eqref{defwpd}, in a domain that is either periodic or bounded with the Dirichlet boundary conditions. We assume that both phases share the same low-Mach-number scaling, i.e., a single Mach number $\epsilon$ characterizes the flow. Taking the limit as $\eps$ goes to zero in the system \eqref{algebraic2velocityeps}, we get formally the following limit system:
\begin{equation}\label{limitsys2vel}
    \begin{cases}
        \rp^0 = C^{\frac{1}{\gp}},\\
        \rmi^0 = C^{\frac{1}{\gm}},\\
        \div (\ap^0 u_+^0 +\am^0 u_-^0) = 0,\\
        \dt (\ap^0\rp^0 u_+^0) + \div (\ap^0\rp^0 u_+^0 \otimes u_+^0) + \ap^0 \nabla \pi^0 + p_{\mathrm{int}}^0 \nabla \ap^0 = 0,\\
        \dt(\am^0 \rmi^0 u_-^0) + \div (\am^0 \rmi^0 u_-^0 \otimes  u_-^0) + \am^0 \nabla \pi^0 + p_{\mathrm{int}}^0 \nabla \am^0 = 0,\\
        \dt \alpha_\pm^0 + \div (\alpha_\pm^0 u_\pm^0) =0,\\
        \ap^0 + \am^0 = 1.\\
    \end{cases}
\end{equation}
with $\pi^0$ the limit pressure.
\end{proposition}

\begin{remark}
In the limit system, both phases are homogeneous ($\rho_\pm$ are constants), but the divergence free condition is now on the averaged velocity $u=\ap u_+ + \am u_-$.
\end{remark}

\begin{remark}
The system \eqref{limitsys2vel} has a similar form than the two-layers shallow-water system between rigid lids when we reformulate using $U=\ap u_+ + \am u_-$ with $\div U =0$, $h_++h_-=\text{cst}$, $p_{\mathrm{int}}=g$ the gravity and with a term $\text{cst} \nabla \ap$ in the momentum equation of phase $+$ (see \cite{bresch2011well}).
\end{remark}

\begin{proof}[Formal proof]
The formal proof is very similar to the one in Part~\ref{formalp} for the two-phase system with one velocity and a pressure algebraic closure. At order $\eps^{-2}$ the momentum equations yield
\begin{align*}
    \left\{
        \begin{array}{l}
              \nabla p_+(\rho_+^0) = 0, \\
              \nabla p_-(\rho_-^0) = 0.
        \end{array}
    \right.
\end{align*}

Following the same development than in the proof of Part~\ref{formalp}, because we also have $p_+(\rp^0)=p_-(\rmi^0)$ from the equality of pressures, we obtain
\begin{align*}
    \rho_+^0 = C_0^{1/\gamma_+},\\
    \rho_-^0 = C_0^{1/\gamma_-}.
\end{align*}

But now the mass equations at order 1 read
\begin{align*}
    \dt \alpha_\pm^0 + \div (\alpha_\pm^0 u_\pm^0) = 0.
\end{align*}

Adding these two equations and using $\ap^0+\am^0=1$ we get
\begin{align*}
    \div (\ap^0 u_+^0 + \am^0 u_-^0) = 0.
\end{align*}

From the algebraic closure at order $\eps^2$, we know that $p_+'(\rp^0) \rp^1 = p_-'(\rmi^0) \rmi^1$. We set $\pi^0 := p_+'(\rp^0) \rp^1$. Therefore, the momentum equations at order 1 read
\begin{align*}
    \rp^0 \ap^0 ( \dt u_+^0 + u_+^0 \cdot\nabla u_+^0) + \ap^0 \nabla \pi^0 + p_{\mathrm{int}}^0 \nabla \ap^0 = 0,\\
    \rmi^0 \am^0 (\dt u_-^0 + u_-^0 \cdot\nabla u_-^0) + \am^0 \nabla \pi^0 + p_{\mathrm{int}}^0 \nabla \am^0 = 0.
\end{align*}
Note that the terms $\alpha_\pm^1 \nabla p_\pm(\rho_\pm^0)$ are cancelled because we have $\nabla p_\pm(\rho_\pm^0)=0$ thanks to the momentum equations at order $\eps^{-1}$.
\end{proof}

\subsection{The non-isentropic case}

\subsubsection{With algebraic closure}

We consider a non-isentropic two-phase flow system, where each phase has its own entropy $S_\pm$. As in the previous section, we neglect the effects of viscosity and surface tension and consider only an interfacial pressure term $p_{\mathrm{int}} \nabla \ap$. We consider the following system with an algebraic closure for the pressure:
\begin{align*}
    \left\{
        \begin{array}{l}
            \ap + \am = 1,\\
            \dt (\ap \rp) + \div(\ap \rp u_+) = 0,\\closure
            \dt (\am \rmi) + \div(\ap \rmi u_-) = 0,\\
            \dt (\ap \rp u_+) + \div (\ap \rp u_+ \otimes u_+) + \ap \nabla p_+(\rp, S_+) + p_{\mathrm{int}} \nabla \ap = 0,\\
            \dt (\am \rmi u_-) + \div (\am \rmi u_- \otimes u_-) + \am \nabla p_-(\rmi, S_-) + p_{\mathrm{int}} \nabla \am = 0,\\
            \dt S_\pm + u_\pm \cdot \nabla S_\pm = 0,\\
            p_+(\rp,S_+)=p_-(\rmi,S_-).
        \end{array}
    \right.
\end{align*}
For simplicity, we consider $p_\pm(\rho_\pm,S_\pm) = \rho_\pm^{\gamma_\pm}e^{S_\pm}$ for the state pressure laws. To the best of our knowledge, there are no works on the global existence of weak solutions for this system.

\begin{remark}
We refer the interested reader to \cite{ndjinga2007influence} for a discussion on the hyperbolicity of this sytem, in particular depending on the term $p_{\mathrm{int},i}$.\newline
\end{remark}

We consider the following adimensional system:
\begin{align}
    \left\{
        \begin{array}{l}\label{PDE2veleps}
            \ap^\eps + \am^\eps = 1,\\
            \dt (\ap^\eps \rp) + \div(\ap^\eps \rp^\eps u_+^\eps) = 0,\\
            \dt (\am^\eps \rmi^\eps) + \div(\ap^\eps \rmi^\eps u_-^\eps) = 0,\\
            \dt (\ap^\eps \rp^\eps u_+^\eps) + \div (\ap^\eps \rp^\eps u_+^\eps \otimes u_+^\eps) + \frac{1}{\eps^2} \ap^\eps \nabla p_+(\rp^\eps, S_+^\eps) + p_{\mathrm{int}} \nabla \ap^\eps = 0,\\
            \dt (\am^\eps \rmi^\eps u_-^\eps) + \div (\am^\eps \rmi^\eps u_-^\eps \otimes u_-^\eps) + \frac{1}{\eps^2} \am^\eps \nabla p_-(\rmi^\eps, S_-^\eps) + p_{\mathrm{int}} \nabla \am^\eps = 0,\\
            \dt S_\pm^\eps + u_\pm^\eps \cdot \nabla S_\pm^\eps = 0,\\
            p_+(\rp^\eps, S_+^\eps) = p_-(\rmi^\eps, S_-^\eps).
        \end{array}
    \right.
\end{align}

In this subsection we are going to give a formal proof of the following proposition.

\begin{proposition}
We consider a smooth enough solution of system \eqref{PDE2veleps} associated with well-prepared initial data in the sense of Definition \eqref{defwpd}, in a domain that is either periodic or bounded with the Dirichlet boundary conditions. We assume that both phases share the same low-Mach-number scaling, i.e., a single Mach number $\epsilon$ characterizes the flow. Taking the limit as $\eps$ goes to zero in the system \eqref{PDE2veleps}, we get formally the following limit system:
\begin{equation*}
    \begin{cases}
        \rp^0 = C^{\frac{1}{\gp}}e^{-S_+^0/\gp},\\
        \rmi^0 = C^{\frac{1}{\gm}}e^{-S_-^0/\gm},\\
        \div (\ap^0 u_+^0 +\am^0 u_-^0) = 0,\\
        \dt (\ap^0 \rp^0 u_+^0) + \div (\ap^0 \rp^0 u_+^0 \otimes u_+^0) + \ap^0 \nabla \pi^0 + p_{\mathrm{int}}^0 \nabla \ap^0 = 0,\\
        \dt (\am^0 \rmi^0 u_-^0) + \div (\am^0 \rmi^0 u_-^0 \otimes u_-^0) + \am^0 \nabla \pi^0 + p_{\mathrm{int}}^0 \nabla \am^0 = 0,\\
        \dt S_\pm^0 + u_\pm^0 \cdot \nabla S_\pm^0=0,\\
        \dt \alpha_+^0 + u_+^0 \cdot \nabla \alpha_+^0 = 0,\\
        \ap^0 + \am^0 = 1,\\
        \dt \alpha_+^0 + \div (\alpha_+^0 u_+^0) = 0.
    \end{cases}
\end{equation*}
with $\pi^0$ the limit pressure (of both phases).
\end{proposition}

\begin{proof}[Formal proof]
This formal proof is very similar to the one in the non-isentropic two-phase one-velocity case and we use the same expansion for the pressure \eqref{expp}. We obtain at order 1:
\begin{align*}
    \ap^0 + \am^0 = 1,\\
    \dt S_\pm^0 + u_\pm^0 \cdot \nabla S_\pm^0 = 0.
\end{align*}
and the algebraic closure reads
\begin{align*}
    p_+ (\rp^0,S_+^0) = p_-(\rmi^0,S_-^0).
\end{align*}

At order $\eps^{-2}$ the momentum equations yield
\begin{align*}
    \left\{
        \begin{array}{l}
              \nabla p_+(\rho_+^0) = 0, \\
              \nabla p_-(\rho_-^0) = 0.
        \end{array}
    \right.
\end{align*}

Thus, we can set $C(t)$ a constant in space such that
\begin{align*}
     p_+ (\rp^0,S_+^0) = p_-(\rmi^0,S_-^0) = C(t),
\end{align*}

Adding the two mass equations using the entropy equations and the state laws, we now obtain
\begin{align}
    \div (\ap^0 u_+^0 + \am^0 u_-^0) + \Big( \frac{\ap^0}{\gp} + \frac{\am^0}{\gm} \Big) \frac{C'(t)}{C(t)}=0.\label{divuI2}
\end{align}

Now we integrate in space, and because of our assumption on $\Omega$, we get
\begin{align*}
    \frac{C'(t)}{C(t)}\int_\Omega \Big( \frac{\ap^0}{\gp} + \frac{\am^0}{\gm} \Big) d x   =0.
\end{align*}

Thus $C'(t) = 0$ and $C$ is a constant in time and space. Coming back to \eqref{divuI2} we obtain
\begin{align*}
    \div (\ap^0 u_+^0 + \am^0 u_-^0) = 0.
\end{align*}

The mass equations at order $1$ now read
\begin{align*}
    \dt \alpha_\pm^0 + \div (\alpha_\pm^0 u_\pm^0)=0.
\end{align*}

Setting $\pi_i^0 = p_i(\rho_i^0) \Big( \gamma_i \frac{\rho_i^1}{\rho_i^0} + S_i^1 \Big)$ for $i=+,-$, we get from the order $\eps^2$ of the algebraic closure:
\begin{align*}
    \partial_\rho p_+(\rp^0,S_+^0)\rp^1 + \partial_S p_+ S_+^1 &= (\partial_\rho p_- \rmi^1 + \partial_S p_- S_-^1),\\
    \pi_+^0 &= \pi_-^0 := \pi^0. 
\end{align*}

Finally, the order 1 momentum equations read
\begin{align*}
    \rp^0 \ap^0 ( \dt u_+^0 + u_+^0 \nabla \cdot u_+^0) + \ap^0 \nabla \pi^0 + p_{\mathrm{int}}^0 \nabla \ap^0 = 0,\\
    \rmi^0 \am^0 (\dt u_-^0 + u_-^0 \nabla \cdot u_-^0) + \am^0 \nabla \pi^0 + p_{\mathrm{int}}^0 \nabla \am^0 = 0.
\end{align*}

\end{proof}

\subsubsection{With PDE closure}

We now consider a non-isentropic two-phase flow system with a PDE closure for the pressure:
\begin{align*}
    \left\{
        \begin{array}{l}
            \ap + \am = 1,\\
            \dt (\ap \rp) + \div(\ap \rp u_+) = 0,\\
            \dt (\am \rmi) + \div(\ap \rmi u_-) = 0,\\
            \dt (\ap \rp u_+) + \div (\ap \rp u_+ \otimes u_+) + \ap \nabla p_+(\rp, S_+) + p_{\mathrm{int}} \nabla \ap = 0,\\
            \dt (\am \rmi u_-) + \div (\am \rmi u_- \otimes u_-) + \am \nabla p_-(\rmi, S_-) + p_{\mathrm{int}} \nabla \am = 0,\\
            \dt S_\pm + u_\pm \cdot \nabla S_\pm = 0,\\
            \dt \ap + u_{\mathrm{int}} \cdot \nabla \ap = \ap \am (p_+(\rp, S_+)-p_-(\rmi, S_-)),
        \end{array}
    \right.
\end{align*}
with $u_{\mathrm{int}}$ the interface velocity that needs to be specified according to the physical context. We assume here that $u_{\mathrm{int}} = u_+$ (see \cite{narbona2024two}). 

\begin{remark}
Note that the choice of the terms $u_{\mathrm{int}}$ and $p_{\mathrm{int}}$ is an important choice for the model because it ensures the hyperbolicity of the system. We refer the interested reader to \cite{coquel2002closure} or \cite{gallouet2010hyperbolic}  among others for more information on this subject.\newline
\end{remark}

We consider the following adimensional system:
\begin{align}\label{noniPDE}
    \left\{
        \begin{array}{l}
            \ap^\eps + \am^\eps = 1,\\
            \dt (\ap^\eps \rp) + \div(\ap^\eps \rp^\eps u_+^\eps) = 0,\\
            \dt (\am^\eps \rmi^\eps) + \div(\ap^\eps \rmi^\eps u_-^\eps) = 0,\\
            \dt (\ap^\eps \rp^\eps u_+^\eps) + \div (\ap^\eps \rp^\eps u_+^\eps \otimes u_+^\eps) + \frac{1}{\eps^2} \ap^\eps \nabla p_+(\rp^\eps, S_+^\eps) + p_{\mathrm{int}} \nabla \ap^\eps = 0,\\
            \dt (\am^\eps \rmi^\eps u_-^\eps) + \div (\am^\eps \rmi^\eps u_-^\eps \otimes u_-^\eps) + \frac{1}{\eps^2} \am^\eps \nabla p_-(\rmi^\eps, S_-^\eps) + p_{\mathrm{int}} \nabla \am^\eps = 0,\\
            \dt S_\pm^\eps + u_\pm^\eps \cdot \nabla S_\pm^\eps = 0,\\
            \dt \ap^\eps + u_+^\eps \cdot \nabla \ap^\eps = \ap^\eps \am^\eps \frac{1}{\eps^2} \Big(p_+(\rp^\eps, S_+^\eps)-p_-(\rmi^\eps, S_-^\eps)\Big).
        \end{array}
    \right.
\end{align}

In this subsection we are going to give a formal proof of the following proposition.

\begin{proposition}
We consider a smooth enough solution of system \eqref{noniPDE} associated with well-prepared initial data in the sense of Definition \eqref{defwpd}, in a domain that is either periodic or bounded with the Dirichlet boundary conditions. We assume that both phases share the same low-Mach-number scaling, i.e., a single Mach number $\epsilon$ characterizes the flow. Taking the limit as $\eps$ goes to zero in the system \eqref{noniPDE}, we get formally the following limit system:
\begin{equation*}
    \begin{cases}
        \rp^0 = C^{\frac{1}{\gp}}e^{-S_+^0/\gp},\\
        \rmi^0 = C^{\frac{1}{\gm}}e^{-S_-^0/\gm},\\
        \div (\ap^0 u_+^0 +\am^0 u_-^0) = 0,\\
        \dt (\ap^0 \rp^0 u_+^0) + \div (\ap^0 \rp^0 u_+^0 \otimes u_+^0) + \ap^0 \nabla \pi_+^0 + p_{\mathrm{int}}^0 \nabla \ap^0 = 0,\\
        \dt (\am^0 \rmi^0 u_-^0) + \div (\am^0 \rmi^0 u_-^0 \otimes u_-^0) + \am^0 \nabla \pi_-^0 + p_{\mathrm{int}}^0 \nabla \am^0 = 0,\\
        \rp^0 \ap^0 ( \dt u_+^0 + u_+^0 \nabla \cdot u_+^0) + \ap^0 \nabla \pi^0 + p_{\mathrm{int}}^0,\\
        \dt S_\pm^0 + u_\pm^0 \cdot \nabla S_\pm^0=0,\\
        \ap^0 \div(u_+^0) = \ap^0 \am^0 (\pi_-^0 - \pi_+^0 ),\\
        \ap^0 + \am^0 = 1,\\
        \dt \alpha_+^0 + \div (\alpha_+^0 u_+^0) = 0.
    \end{cases}
\end{equation*}
with $\pi_i^0$ the limit pressure of phase $i$, $i=+,-$.
\end{proposition}

\begin{remark}
    Like in \cite{narbona2024two}, we do not get the equality of the limit pressures, but a relation between them involving $\div u_+^0$, in the form: 
    $$\ap^0 \div(u_+^0) = \ap^0 \am^0 (\pi_-^0 - \pi_+^0 ).$$
    This allows to consider only one of the pressures as a variable, the other being given by this constraint.
\end{remark}

\begin{proof}[Formal proof]
This formal proof is very similar to the one in the non-isentropic two-phase two-velocity case with algebraic closure. In the same way, we obtain
\begin{align*}
     p_+ (\rp^0,S_+^0) = p_-(\rmi^0,S_-^0) = C,\\
     \div (\ap^0 u_+^0 + \am^0 u_-^0) = 0.
\end{align*}

Note that the PDE closure at order $\eps^{-2}$ reads as the algebraic closure
\begin{align*}
    p_+ (\rp^0,S_+^0) = p_-(\rmi^0,S_-^0).
\end{align*}

The mass equations at order $1$ now read
\begin{align*}
    \dt \alpha_\pm^0 + \div (\alpha_\pm^0 u_\pm^0)=0.
\end{align*}

Setting $\pi_i^0 = p_i(\rho_i^0,S_i^0) \Big( \gamma_i \frac{\rho_i^1}{\rho_i^0} + S_i^1 \Big)$ for $i=+,-$, we get from the order $1$ of the PDE closure:
\begin{align*}
    \dt \ap^0 + u_+^0 \cdot \nabla \ap^0 &= \ap^0 \am^0 \Big( \partial_\rho p_+(\rp^0,S_+^0)\rp^1 + \partial_S p_+ S_+^1 - (\partial_\rho p_- \rmi^1 + \partial_S p_- S_-^1) \Big),\\
    \dt \ap^0 + u_+^0 \cdot \div \ap^0 &= \ap^0 \am^0 (\pi_+^0 - \pi_-^0. ).
\end{align*}

Using the mass equation to replace $\dt \ap^0$, we get
\begin{align*}
    -\ap^0 \div(u_+^0) = \ap^0 \am^0 (\pi_+^0 - \pi_-^0 ).
\end{align*}

Finally, the order 1 momentum equations read
\begin{align*}
    \rp^0 \ap^0 ( \dt u_+^0 + u_+^0 \nabla \cdot u_+^0) + \ap^0 \nabla \pi_+^0 + p_{\mathrm{int}}^0 \nabla \ap^0 = 0,\\
    \rmi^0 \am^0 (\dt u_-^0 + u_-^0 \nabla \cdot u_-^0) + \am^0 \nabla \pi_-^0 + p_{\mathrm{int}}^0 \nabla \am^0 = 0.
\end{align*}

\end{proof}

\section{Conclusion and open problems}\label{Opb}

The formal analysis of the low-Mach-number limit for two-phase flows has made it possible to characterize the resulting limit systems according to the increasing complexity of the models considered. For a single-velocity system, the absence of entropy leads to non-homogeneous incompressible limit systems, while the introduction of entropy alters the structure of the equations and the nature of the coupling between phases. In the two-velocity case, results differ depending on whether entropy is considered and on the type of closure: for isentropic systems, algebraic and PDE closures yield distinct limit systems. These findings highlight the importance of variable selection and closure assumptions on the form and stability of reduced models, and call for further investigations. We now discuss several open problems, which naturally arise from this analysis and represent key research directions.

\textbf{A non-isentropic case with more complicated thermodynamics.} In this paper, we consider a system of equations for non-isentropic flows with the simplest entropy equations. This allows for an easier formal analysis of the low Mach number limit. One could improve the physical basis of the system by adding thermodynamics, using temperatures for each phase to get additional energy equations. This is done by \cite{varsakelis2011low} that generalize the low-Mach-number asymptotics for two-phase flows of granular materials. After non-dimensionalizing the system with respect to a reference thermodynamic state, the Mach number appears in the energy equations. Considering a PDE closure for the pressure that is an equation on the volume fraction that depends on both pressure difference and phasial temperatures difference, they get in the limit equality of both pressures and temperatures. Then only one energy equation is needed in the system, and the other can be used as an equation for the volume fraction. This is particularly advantageous in a numerical point of view.\\

\textbf{A different Mach number for each phase.} In this paper we have always considered the same Mach number for both phases. To get closer to the physics, one would like to consider the low-Mach-number limit of only one Mach number, the one of the phase with the smallest speed of sound. The asymptotic calculation is done by \cite{narbona2024two}, where the authors use a system of equations governing two compressible gas phases, following respectively a perfect gas law and a stiff gas pressure state law, and perform a low-Mach-number limit with respect to the second phase. They obtain in the limit a constant density for the second phase, but not anymore equality of pressures or of temperatures nor the free divergence of the velocity of the second phase.\\

\textbf{General initial data in more general domains.} In this paper, we only presented some results for the low Mach number limit for multiphase flows in the case of well-prepared data, in a bounded domain. As far as we know, nothing has been done in the case of more general initial data in other types of domains.\\

\textit{Acknowledgement.} The author would like to thank Didier Bresch and Christophe Lacave for their many tips on writing this text.
She also gratefully acknowledges the University Savoie Mont Blanc and the Agence Nationale pour la Recherche (ANR) for her PhD position fellowship linked to the CPJ Anamod managed by Christophe Lacave. She also wants to thank for some support of the ANR under France 2030 bearing the reference ANR-23-EXMA-004 (Complexflows project) related to the PEPR Maths-ViVEs. The author wants also to gratefully acknowledge the partial support by the ANR grant ANR-23-CE40-0014-01 (ANR Bourgeons).

%\bibliographystyle{apalike}
%\bibliography{bibli}

\adrese
\end{document}